\documentclass[11pt,twoside,reqno]{amsart}

\usepackage{amsmath}
\usepackage{amsthm}
\usepackage{amsfonts,amssymb}
\usepackage{latexsym}
\usepackage{mathrsfs}
\usepackage[all]{xy}
\usepackage{url}
\usepackage[pdftex]{graphicx}
\usepackage{float}
\usepackage{hyperref}

\theoremstyle{plain}
\newtheorem{lema}{Lemma}[section]
\newtheorem{prop}[lema]{Proposition}
\newtheorem{teo}[lema]{Theorem}
\newtheorem*{teoppal}{Theorem \ref{CasoAbeliano}}

\newtheorem{preg}[lema]{Question}
\newtheorem{coro}[lema]{Corollary}
\theoremstyle{remark}

\theoremstyle{definition}
\newtheorem{defi}[lema]{Definition}

\newcommand{\deficiency}{\mathrm{def}}
\newcommand{\Fix}{\mathrm{Fix}}

\newcommand{\Z}{\mathbb{Z}}
\newcommand{\Q}{\mathbb{Q}}

\newcommand{\PP}{\mathcal{P}}

\newcommand{\RR}{\mathcal{R}}
\newcommand{\TT}{\mathcal{T}}

\newcommand{\wt}{\widetilde}
\def\Re{\mathrm{Re}}
\def\Im{\mathrm{Im}}

\begin{document}

\title[On a question of R.H. Bing]{On a question of R.H. Bing concerning the fixed point property for two-dimensional polyhedra}

\author[J.A. Barmak]{Jonathan Ariel Barmak $^{\dagger}$}
\author[I. Sadofschi Costa]{Iv\'an Sadofschi Costa}

\thanks{$^{\dagger}$ Researcher of CONICET. Partially supported by grants ANPCyT PICT-2011-0812, CONICET PIP 112-201101-00746 and UBACyT 20020130100369.}

\address{Departamento de Matematica\\
FCEyN-Universidad de Buenos Aires\\
Buenos Aires, Argentina}

\email{jbarmak@dm.uba.ar}
\email{ivansadofschi@gmail.com}

\begin{abstract}
In 1969 R.H. Bing asked the following question: Is there a compact two-dimensional polyhedron with the fixed point property which has even Euler characteristic? In this paper we prove that there are no spaces with these properties and abelian fundamental group. We also show that the  fundamental group of an example cannot have trivial Schur multiplier. 
\end{abstract}

\subjclass[2010]{55M20, 57M20, 55P15, 57M05}

\keywords{Fixed point property, two-dimensional complexes, Nielsen fixed point theory, homotopy classification.}

\maketitle

\section{Introduction}

In his famous article ``The elusive fixed point property" of 1969 \cite{Bing}, R.H. Bing stated twelve questions which were a motivation for the development of several methods in the fixed point theory of polyhedra and continua. In the last 45 years, eight of these questions have been answered while four of them remain still open \cite{Hagopian}. In this paper we make an advance on Bing's Question 1.

Recall that a space $X$ is said to have the \textit{fixed point property} if every continuous self-map $f:X\to X$ has a fixed point. The fixed point property is clearly a topological invariant, but W. Lopez showed in \cite{Lopez} that it is not a homotopy invariant in the category of compact polyhedra (although it is \textit{almost} a homotopy invariant, see Theorem \ref{JiangFPPInvariante} below). In order to provide an example Lopez constructed an eight-dimensional polyhedron with the fixed point property and even Euler characteristic. This example motivated the following question.

\begin{preg}[Bing's Question 1]\label{Bing1} Is there a compact two-dimensional polyhedron with the fixed point property which has even Euler characteristic?  
\end{preg}

It is not hard to find examples of two-dimensional complexes with the fixed point property. For any finite group $G$ with deficiency equal to $0$, there exists a $2$-complex $X$ with $\pi_1(X)=G$ and $\wt{H}_*(X;\Q)=0$. By the Lefschetz fixed point theorem, any compact two-dimensional polyhedron $X$ with trivial rational homology has the fixed point property. However, it is unknown whether the converse of the latter statement holds. Therefore
we will also consider the following variation of Bing's question:

\begin{preg}\label{Bing2} Is there a compact two-dimensional polyhedron with the fixed point property such that $\wt{H}_*(X;\Q)\neq 0$?
\end{preg}

Of course, an affirmative answer to Question \ref{Bing1} implies an affirmative answer to Question \ref{Bing2}.

It is well-known that a polyhedron $X$ with $H_1(X;\mathbb{Q})\neq 0$ lacks the fixed point property since $S^1$ is a retract of any such space. Therefore, for a compact $2$-complex with the fixed point property it is equivalent to saying that $\wt{H}_*(X;\mathbb{Q})\neq 0$, that $\chi (X)> 1$ or that $H_2(X)\neq 0$. 

A higher dimensional analogue to Question \ref{Bing2} has been settled by Waggoner \cite{Waggoner} for dimension $n\geq 4$ and later extended to dimension $3$ by Jiang \cite{Jiang}.

\begin{teo}[Waggoner, Jiang] \label{Waggoner} If $X$ is a compact $(n-2)$-connected polyhedron of dimension $n>2$ and $\wt{H}_*(X;\Q)\neq 0$, then $X$ does not have the fixed point property.
\end{teo}

In this article a compact two-dimensional polyhedron will be called a \textit{Bing space} if it has the fixed point property and $\wt{H}_*(X;\mathbb{Q})\neq 0$. The main result of this paper is the following

\begin{teoppal}
There are no Bing spaces with abelian fundamental group.
\end{teoppal}

In particular, no space of those considered by Question 1.1 can have abelian fundamental group. The ideas involved in the proof include the relationship between $2$-complexes and group presentations, the homotopy classification of $2$-complexes with finite abelian fundamental group, Nielsen fixed point theory, elementary results on two-dimensional homotopy theory and obstruction theory.

The strategy in the proof will also be used to prove that there are no Bing spaces with fundamental group isomorphic to the alternating groups $A_4, A_5$, the symmetric group $S_4$ or to any dihedral group.

The first part of the paper is devoted to proving that there are no Bing spaces whose fundamental group has trivial Schur multiplier. 

\section{Preliminaries}
\label{Prelim}

We will denote by $H_*(X)$ the homology of $X$ with integral coefficients. A Bing space must necessarily be path-connected, therefore the basepoint will be omitted in the notation for fundamental groups. By the comments above, if $X$ is a Bing space, $H_1(X)$ must be a torsion group, since otherwise $S^1$ would be a retract of $X$, and the fixed point property is preserved by retracts. It is easy to see that $X\vee Y$ has the fixed point property if and only if both $X$ and $Y$ have the fixed point property.

\begin{defi}
Let $X$ be a connected polyhedron. We say that $x\in X$ is a \textit{local separating point} if there is a connected open neighborhood $U$ of $x$ such that $U-\{x\}$ is not connected. We say that $x\in X$ is a \textit{global separating point} if $X-\{x\}$ is not connected. In particular, any global separating point is a local separating point.
\end{defi}

Note that for a given triangulation every local separating point is a vertex or lies in a maximal $1$-simplex.   

The following result was proved by Jiang in \cite[Theorem 7.1]{Jiang}.
\begin{teo}[Jiang]\label{JiangFPPInvariante}
In the category of compact connected polyhedra without global separating points, the fixed point property is a homotopy type invariant.

Moreover, if $X\simeq Y$ are compact connected polyhedra such that $Y$ lacks the fixed point property and $X$ does not have global separating points, then $X$ lacks the fixed point property.
\end{teo}
The second part of the result above does not appear in Jiang's original formulation. It can be deduced from its proof or, alternatively, from the first part by replacing $Y$ by $Y\times I$.

Recall that a presentation $\PP=\langle a_1,\ldots, a_n\mid r_1,\ldots , r_k \rangle$ of a group $G$ has associated a CW-complex $K_{\PP}$, called the \textit{standard complex} of $\PP$, with one $0$-cell, a $1$-cell for each generator $a_i$ and a $2$-cell for every relator $r_j$. The fundamental group of $K_{\PP}$ is $G$. The second barycentric subdivision of this complex is a triangulation of $K_{\PP}$, so any standard complex is in fact a polyhedron. Conversely, every compact connected $2$-complex $X$ is homotopy equivalent to the standard complex of a presentation $\PP$ of $\pi_1(X)$.

\begin{defi} The \textit{deficiency} of a presentation $\PP=\langle a_1,\ldots,a_n \mid r_1,\ldots , r_k \rangle$ is defined by $\deficiency(\PP)=k-n$. Therefore, $\chi(K_{\PP})=\deficiency(\PP)+1$. Given a finitely presented group $G$, its \textit{deficiency} $\deficiency(G)$ is the minimum possible deficiency of a presentation of $G$. Then, for any compact connected $2$-complex $X$ with fundamental group $G$, $\chi (X)\ge \deficiency(G)+1$. We say that $X$ has \textit{minimum Euler characteristic} if $\chi(X)=\deficiency(G)+1$.
\end{defi}

For every connected CW-complex $X$ there is a short exact sequence (\cite[Chapter II, Theorem 5.2]{BrownCohomology})
$$0\to \Sigma_2(X)\to H_2(X)\to H_2(G)\to 0.$$
Here $\Sigma_2(X)$ stands for the subgroup of spherical cycles in $H_2(X)$, that is the image of the Hurewicz map $h:\pi_2(X)\to H_2(X)$, and $G=\pi_1(X)$. The second homology group $H_2(G)$ of $G$ is called the \textit{Schur multiplier} of $G$. When $G$ is a finitely presented group, the sequence above provides a lower bound for the deficiency of $G$. If in addition $H_1(G)$ is a torsion group, we have
$$\deficiency(G)\geq \text{number of invariant factors of } H_2(G).$$
If the bound above is sharp we say that $G$ is an \textit{efficient} group.

\section{Primitive spherical elements and Waggoner's Theorem}

A strategy for proving that a space $X$ without global separating points lacks the fixed point property is to show that there exists a space $Y\simeq X$ that has $S^n$ as a retract. Waggoner used this idea to prove Theorem \ref{Waggoner}, which involves only simply-connected spaces. In this section we show to what extent we can apply this strategy in the case of an arbitrary $2$-complex. The exact sequence of spherical elements will be used instead of the Hurewicz Theorem to characterize precisely the situations in which this idea can be applied. In particular we will deduce that there are no Bing spaces with fundamental group $G$ if $H_2(G)=0$. We will also show that if the fundamental group $G$ of a Bing space $X$ is freely indecomposable, then $X$ must have minimum Euler characteristic and $G$ must be efficient.

The following result appears essentially in \cite{Waggoner}. We exhibit here a shorter proof.

\begin{lema}[Waggoner]\label{LemaWaggoner} Let $(X,S^n)$ be a CW-pair with $\dim(X)\leq n+1$, $n\geq 1$. If $i_*:H_{n}(S^n)\to H_n(X)$ is a split monomorphism, then $S^n$ is a retract of $X$.
\begin{proof}
By the naturality of the short exact sequence in the universal coefficient theorem $i^*:H^n(X;\pi_n(S^n))\to H^n(S^n;\pi_n(S^n))$ is surjective. Then, the connecting homomorphism $\delta: H^n(S^n; \pi _n (S^n))\to H^{n+1}(X,S^n; \pi _n (S^n))$ is trivial. By obstruction theory \cite[Theorem 8.4.1]{Spanier}, $i:S^n\hookrightarrow X$ can be extended to $X$.
\end{proof}
\end{lema}

Let $F$ be a free abelian group. We say that $a\in F$ is \textit{primitive} in $F$ if the homomorphism $\Z\to F$ defined by $1\mapsto a$ is a split monomorphism. This is equivalent to saying that $a$ can be extended to a basis of $F$.

The proof of the following lemma is an easy application of the Smith normal form.
\begin{lema}\label{LemaSmith}
Consider an exact sequence $0 \to S \to \Z^k \to A \to 0$ of abelian groups. Then the number of invariant factors of $A$ is strictly smaller than $k$ if and only if there exists $a\in S$ primitive in $\Z^k$.
\end{lema}

\begin{prop}\label{CuandoS2EsRetracto}
Let $X$ be a compact connected $2$-dimensional polyhedron. The following are equivalent:

(i) $X$ is homotopy equivalent to a polyhedron $Y$ having $S^2$ as a retract.

(ii) There exists $a\in \Sigma_2(X)$ primitive in $H_2(X)$.

(iii) The number of invariant factors of $H_2(\pi_1(X))$ is strictly smaller than the rank of $H_2(X)$.
\begin{proof}
If $S^2$ is a retract of a space $Y$ homotopy equivalent to $X$, there are maps $f:S^2\to X$ and $g:X\to S^2$ such that $gf\simeq 1_{S^2}$. Then $h([f])\in \Sigma_2 (X)$ is a primitive element of $H_2(X)$. Conversely, if $\Sigma _2(X)$ contains a primitive element $h([f])$ of $H_2(X)$, then $f_*:H_2(S^2)\to H_2 (X)$ is a split monomorphism. We may assume that $f$ is a simplicial map for some triangulations of $S^2$ and $X$, so the mapping cylinder $Y=M(f)$ is a three-dimensional polyhedron. The canonical inclusion $i:S^2\hookrightarrow Y$ induces a split monomorphism in $H_2$ and Lemma \ref{LemaWaggoner} says then that $S^2$ is a retract of $Y$.

The equivalence between (ii) and (iii) is immediate from Lemma \ref{LemaSmith}.
\end{proof}
\end{prop}

\begin{teo}\label{H2Trivial} If $H_2(G)=0$ there are no Bing spaces with fundamental group $G$.
\begin{proof}
Suppose $X$ is a Bing space with $\pi_1(X)=G$. Note that $X=X_1\vee \ldots \vee X_m$, where each $X_i$ is a polyhedron without global separating points or a $1$-simplex (the basepoints of the wedges may not be the same). For some $i$, we must have $H_2(X_i)\neq 0$. Since $X_i$ is a retract of $X$, it must have the fixed point property. On the other hand, $\pi _1(X_i)$ is a free factor of $G$, and then $H_2(G)=0$ implies $H_2(\pi_1 (X_i))=0$ (\cite[Corollary 6.2.10]{Weibel}). By Proposition \ref{CuandoS2EsRetracto} and Theorem \ref{JiangFPPInvariante}, $X_i$ lacks the fixed point property, a contradiction.
\end{proof}
\end{teo}

The particular case $G=0$ in Theorem \ref{H2Trivial} was previously studied by Waggoner in \cite{Waggoner2}. 

\begin{coro}\label{Largo}
There are no Bing spaces with fundamental group isomorphic to the trivial group, cyclic groups, dihedral groups of order $2\pmod 4$, $\mathrm{SL}(n,\mathbb{F}_q)$ (for $(n,q)\neq (2,4)$, $(2,9)$, $(3,2)$, $(3,4)$, $(4,2)$), deficiency-zero groups (e.g. the quaternion group), groups of square-free order (more generally, any group in which every Sylow subgroup has trivial Schur multiplier), $13$ of the $26$ sporadic simple groups and many infinite families of finite simple groups of Lie type.
\begin{proof}
All these groups have trivial Schur multiplier. For cyclic groups, dihedral groups and $\mathrm{SL}(n,\mathbb{F}_q)$ this appears in \cite{Weibel}. For deficiency-zero groups it follows from the bound given in section \ref{Prelim}. For groups in which every Sylow subgroup has trivial Schur multiplier, it follows from \cite[Chapter III, Corollary 10.2 and Theorem 10.3]{BrownCohomology}. For the statement about finite simple groups, see \cite[Section 6.1]{GLS}.
\end{proof}
\end{coro}

A group $G$ is said to be \textit{freely indecomposable} if $G\simeq H*K$ implies $H\simeq 1 \text{ or } K\simeq 1$. Finite groups and abelian groups clearly are freely indecomposable.

The following reduction will also be helpful in section \ref{SeccionAbelianos}.

\begin{prop}\label{Reduccion} Let $X$ be a Bing space with freely indecomposable fundamental group $G$. Then there is a Bing space $Y\simeq X$ without global separating points.
\begin{proof}
Fix a triangulation of $X$. If $X$ has a global separating point and is not a $1$-simplex, then $X$ is a wedge of two polyhedra $X_1,X_2$, each with fewer vertices than $X$. By van-Kampen's theorem $G\simeq \pi_1(X_1)*\pi_1(X_2)$ and since $G$ is freely indecomposable, one of these two polyhedra, say $X_2$, is simply-connected. By Theorem \ref{H2Trivial} there are no simply-connected Bing spaces, so $\wt{H}_*(X_2)=0$. Therefore $X_2$ is contractible and then $X_1\simeq X$. By induction there exists a Bing space $Y\simeq X_1$ without global separating points.
\end{proof}
\end{prop}

\begin{prop} If $G$ is freely indecomposable, and $X$ is a Bing space with fundamental group $G$, the rank of $H_2(X)$ must equal the number of invariant factors of $H_2(G)$.
\begin{proof}
By Proposition \ref{Reduccion} we may assume $X$ does not have global separating points. If the rank of $H_2(X)$ is strictly greater than the number of invariant factors of $H_2(G)$, by Proposition \ref{CuandoS2EsRetracto} and Theorem \ref{JiangFPPInvariante} we get a contradiction.
\end{proof}
\end{prop}

We deduce the following
\begin{coro}\label{CasoChiNoMinima}
Let $G$ be a freely indecomposable group. Suppose $X$ is a Bing space with fundamental group $G$. Then $G$ is efficient and $X$ has minimum Euler characteristic.
\end{coro}

Finite abelian groups which are non-cyclic are efficient and have non-trivial Schur multiplier. A different strategy will be developed in the
next section to deal with these cases.

\section{Two-complexes with abelian fundamental group}\label{SeccionAbelianos}

In this section we prove that there are no Bing spaces with abelian fundamental group. We will need some basic concepts from Nielsen theory, namely fixed point class, fixed point index and the Nielsen number. We refer the reader to \cite{JiangLibro}, \cite{Brown} for definitions and basic results. We will also need the following theorem of Jiang:

\begin{teo}[Jiang, {\cite[Main Theorem]{Jiang}}]\label{JiangNIgualM}
Let $X$ be a compact connected polyhedron and $f:X\to X$ be a continuous map. If $X$ does not have local separating points and $X$ is not a $2$-manifold (with or without boundary), then there exists $g\simeq f$ with $\# \Fix(g)=N(f)$.
\end{teo}

The following case will be central in our argument.
\begin{lema}\label{FamiliaSinFPP} Let $\PP=\langle a,b \mid a^m, b^n, [a,b]\rangle$. Then $K_{\PP}$ does not have the fixed point property.
\begin{proof}
By Theorem \ref{JiangNIgualM} it suffices to find a map $f:K_\PP\to K_\PP$ such that $N(f)=0$. Let $T=S^1\times S^1\subseteq \mathbb{C}\times \mathbb{C}$. The complex $K_\PP$ can be identified with the following pushout:
\begin{displaymath}
\xymatrix@C=20pt{ S^1_a\coprod S^1_b \ar@{->}^(.55){ (z^m,1) \coprod (1,z^n) }[rrr] \ar@{->}[d] & & & T \ar@{->}^{i_T}[d]\\
		  D^2_a \coprod D^2_b \ar@{->}_(.55){{i_a\coprod i_b}}[rrr] & & & K_\PP } 
\end{displaymath}

Here $S^1_a, S^1_b, D^2_a,D^2_b\subseteq \mathbb{C}$ denote copies of $1$-dimensional spheres and $2$-dimensional disks.

We define $f_T:T\to K_\PP$ by $$f_T(z,w)=i_T\left(-z,-\overline{w}\right).$$

Now we define $f_a:D_a^2\to K_\PP$ and $f_b:D_b^2\to K_\PP$ by
$$
f_a(z)=
\begin{cases}
i_a(2z) & \text{ if } 0\leq |z|\leq \frac{1}{2} \\
i_T\left(\frac{z^m}{|z|^m} \exp(i\pi (2|z|-1)) , \exp(i\pi (2|z|-1))\right)  &\text{ if } \frac{1}{2}\leq |z| \leq 1
\end{cases}
$$
$$
f_b(z)=
\begin{cases}
i_b(2\overline{z}) & \text{ if } 0\leq |z|\leq \frac{1}{2} \\
i_T\left(\exp(i\pi (2|z|-1)) , \frac{\overline{z}^n}{|z|^n} \exp(i\pi (2|z|-1)) \right)  &\text{ if } \frac{1}{2}\leq |z| \leq 1
\end{cases}
$$
A simple verification shows that $f_T, f_a$ and $f_b$ are well-defined and continuous and that they determine a continuous map $f:K_\PP\to K_\PP$.

It is easy to see that the only fixed points of $f$ are $i_a(0)$ and $i_b(0)$. We will show that the two fixed points are in the same fixed point class. Concretely, we exhibit a path $c$ from $i_a(0)$ to $i_b(0)$ such that $c$ and $f\circ c$ are homotopic. Consider the paths $\gamma_a, \delta_a, \delta_b, \gamma_b:[0,1]\to K_\PP$ defined by
\begin{align*}
\gamma_a(t)& =i_a\left(t/2\right) \\
\delta_a(t)&=i_a\left(1/2 + t/2\right) \\
\delta_b(t)&=i_b\left(1- t/2\right) \\
\gamma_b(t)& =i_b\left(1/2-t/2\right)
\end{align*}
The concatenation $\gamma_a*\delta_a*\delta_b*\gamma_b$ is a well-defined path from $i_a(0)$ to $i_b(0)$. 
In order to prove
$$\gamma_a*\delta_a*\delta_b*\gamma_b \simeq f\circ(\gamma_a*\delta_a*\delta_b*\gamma_b)$$
it suffices to show that
\begin{align}
\gamma_a * \delta_a & =  f \circ \gamma_a \\
\delta_b * \gamma_b & = f  \circ \gamma_b \\
e_{i_T(1,1)} &\simeq (f\circ \delta_a) * (f\circ \delta_b),
\end{align}
where $e_{i_T(1,1)}$ denotes the constant loop at $i_T(1,1)$. Equalities (1) and (2) are clear, (3) follows from
$$(f\circ\delta_a)(t)=i_T\left(\exp(i\pi t), \exp(i\pi t)\right)= (f\circ\delta_b)(1-t).$$

Now we show that the unique fixed point class of $f$ is inessential (i.e. the fixed point indices of $i_a(0)$ and $i_b(0)$ add up to zero). One way to see this is by noting that the fixed point indexes of $i_a(0)$ and $i_b(0)$ are $1$ and $-1$ respectively.
Another way is proving that the Lefschetz number  $\Lambda(f)$ is $0$ and invoking the Lefschetz-Hopf theorem \cite[VII, Proposition 6.6]{Dold}. 
% Another way is proving that the Lefschetz number $\Lambda(f)$ is $0$ and invoking the Lefschetz-Hopf theorem which states $\sum_{x\in \Fix(f)}i(f,x)=\Lambda(f)$.
Then, $N(f)=0$ and by Theorem \ref{JiangNIgualM} there exists $g\simeq f$ without fixed points.
\end{proof}
\end{lema}

If $G$ is any finite group, above the minimum Euler characteristic all $2$-complexes with fundamental group $G$ are homotopy equivalent. That fact along with Theorem \ref{clasificacion} below constitutes the classification of homotopy types of compact $2$-complexes with finite abelian fundamental group. We refer to \cite[Chapter III]{TwoDimensional} and \cite{GutierrezLatiolais} for a detailed exposition on this topic.

\begin{teo}[Browning, {\cite[Chapter III, Theorem 2.11]{TwoDimensional}}]\label{clasificacion}
Let $G$ be a finite abelian group with invariant factors $m_1\mid m_2\mid \ldots \mid m_n$. The number of homotopy types of compact connected $2$-complexes with fundamental group $G$ and minimum Euler characteristic is $\left|\Z_{m_1}^*/\pm (\Z_{m_1}^*)^{n-1}\right|$. Every such complex is homotopy equivalent to the standard complex of a presentation
$$\TT_d=\langle a_1,\ldots,a_n \mid a_1^{m_1}, \ldots , a_n^{m_n}, [a_1^d,a_2], [a_i,a_j], i < j, (i,j)\neq(1,2) \rangle$$
with $(d,m_1)=1$.
\end{teo}

\begin{coro}\label{clasificacion2factores} Let $G$ be a finite abelian group with invariant factors $m_1\mid m_2$ and let $X$ be a compact connected $2$-complex with $\pi_1(X)=G$. If $X$ has minimum Euler characteristic, then $X\simeq K_{\PP}$ where $\PP=\langle a_1,a_2 \mid a_1^{m_1}, a_2^{m_2}, [a_1,a_2]\rangle$.
\end{coro}

The last result we need for the proof of Theorem \ref{CasoAbeliano} is the following

\begin{lema}\label{retracto}Let
$$\TT_d=\langle a_1,\ldots,a_n \mid a_1^{m_1}, \ldots ,a_n^{m_n}, [a_1^d,a_2], [a_i,a_j], i < j, (i,j)\neq(1,2) \rangle$$
with $n\ge 2$ and
$$\RR_d=\langle a_1,a_2 \mid a_1^{m_1},  a_2^{m_2}, [a_1^d,a_2] \rangle .$$
Then $K_{\RR_d}$ is a retract of $K_{\TT_d}$.
\begin{proof}

Clearly $K_{\RR_d}$ is a subcomplex of $K_{\TT_d}$. We will define a cellular retraction $r:K_{\TT_d}\to K_{\RR_d}$. The unique $0$-cell of $K_{\TT_d}$ and the $1$-cells $a_1$, $a_2$ are fixed by $r$. The remaining $1$-cells $a_3,\ldots, a_n$ are mapped to the $0$-cell. In the $2$-skeleton, $r$ fixes the $2$-cells $a_1^{m_1}$, $a_2^{m_2}$ and $[a_1^d,a_2]$, and we must extend $r$ to the remaining $2$-cells. This can be achieved since the composition of $r$ with the attaching maps of those cells is null-homotopic.
\end{proof}
\end{lema}

\begin{teo}
\label{CasoAbeliano}
There are no Bing spaces with abelian fundamental group.
\begin{proof}
Let $X$ be a Bing space with abelian fundamental group $G$. Since $H_1(X)=G$ is a finitely generated torsion group, $G$ is finite abelian. Let $m_1\mid m_2\mid \ldots \mid m_n$ be its invariant factors.

Since $G$ is freely indecomposable, by Proposition \ref{Reduccion} we can assume $X$ does not have global separating points. By Corollary \ref{CasoChiNoMinima} we know that $X$ has minimum Euler characteristic. By Theorem \ref{H2Trivial}, $G$ is not cyclic, so $n\geq 2$. From Theorem \ref{clasificacion}, there is a presentation
$$\TT_d=\langle a_1,\ldots,a_n \mid a_1^{m_1}, \ldots ,a_n^{m_n}, [a_1^d,a_2], [a_i,a_j], i < j, (i,j)\neq(1,2) \rangle$$
with $(d,m_1)=1$ such that $X\simeq K_{\TT_d}$. By Theorem \ref{JiangFPPInvariante}, $K_{\TT_d}$ has the fixed point property. Let
$$\RR_d=\langle a_1,a_2 \mid a_1^{m_1}, a_2^{m_2}, [a_1^d,a_2]\rangle .$$
By Lemma \ref{retracto}, $K_{\RR_d}$ is a retract of $K_{\TT_d}$ so $K_{\RR_d}$ has the fixed point property. Finally consider
$$\RR_1=\langle a_1,a_2 \mid a_1^{m_1}, a_2^{m_2}, [a_1,a_2]\rangle .$$
By Corollary \ref{clasificacion2factores}, $K_{\RR_1}\simeq K_{\RR_d}$, therefore by Theorem \ref{JiangFPPInvariante}, $K_{\RR_1}$ has the fixed point property, contradicting Lemma \ref{FamiliaSinFPP}.
\end{proof}
\end{teo}

The ideas used in the proof of the last result can be applied to other cases. The classification of $2$-complexes has been achieved for a few finite groups, aside finite abelian groups. Our last theorem relies on a result of Hambleton and Kreck. In the proof we will use the following

\begin{lema}\label{LemaIndice}
Let $X$ be a compact polyhedron, $f:X\to X$ a map and $F$ a fixed point class of $f$. Suppose there is a subspace $K\subseteq X$ which is itself a compact polyhedron that satisfies:
\begin{itemize}
	\item $f(K)\subseteq K$.
	\item $K$ deformation retracts to $F$.
	\item $F\subseteq K^\circ$.
	\item $F=K\cap \Fix(f)$.
\end{itemize}
Then the fixed point index of $F$ equals its Euler characteristic, that is $i(f,F)=\chi(F)$.
\begin{proof} Let $U=K^\circ$. The fixed point index of $F$ is given by $i(f,F)=i(f|_U:U\to X)$ (\cite[I, Definition 3.8 and Section 4]{JiangLibro}).
We have $$i(f|_U:U\to X)=i(f|_U:U\to K)=i(f|_K:K\to K)=\Lambda(f|_K)=\Lambda(1_F)=\chi(F).$$ The first equality follows from the definition of the fixed point index (\cite[VII, Proposition 5.10]{Dold}). The second equality follows from \cite[VII, (5.11)]{Dold}, the third from the Lefschetz-Hopf theorem (\cite[VII, Proposition 6.6]{Dold}) and the fourth from the fact that $F\hookrightarrow K$ induces isomorphisms in homology.
\end{proof}
\end{lema}

\begin{teo}
There are no Bing spaces with fundamental group $A_4$, $S_4$, $A_5$ or $D_{2n}$.
\begin{proof}

By \cite[Theorem 2.1]{HambletonKreck} for these groups, the homotopy type of a $2$-complex is determined by the Euler characteristic. Consider the following presentations with deficiency $1$:
\begin{align*}
A_4 &=\langle a,b,c\mid a^2, b^3, c^3, abc\rangle, \\
S_4 &=\langle a,b,c\mid a^2, b^3, c^4, abc\rangle, \\
A_5 &=\langle a,b,c\mid a^2, b^3, c^5, abc\rangle, \\
D_{2n} &=\langle a,b,c\mid a^2, b^2, c^n, abc\rangle .
\end{align*}
We only need to prove that the complexes associated to these presentations lack the fixed point property (we do not need to check whether these presentations have minimum deficiency or not). 

Let $\PP=\langle a,b,c\mid a^l, b^m, c^n, abc\rangle$. Consider the space $X=X(l,m,n)$ obtained by deleting three disjoint disks from $S^2$ and then gluing three $2$-cells on the boundaries of these disks, with attaching maps of degrees $l$, $m$ and $n$ (Figure \ref{figura}). We note that $K_\PP$ is a quotient of $X$ by a contractible subcomplex, therefore $K_\PP\simeq X$. We will show $X$ lacks the fixed point property.

\begin{figure}[H] 
\begin{minipage}{6cm}
\includegraphics[scale=0.7]{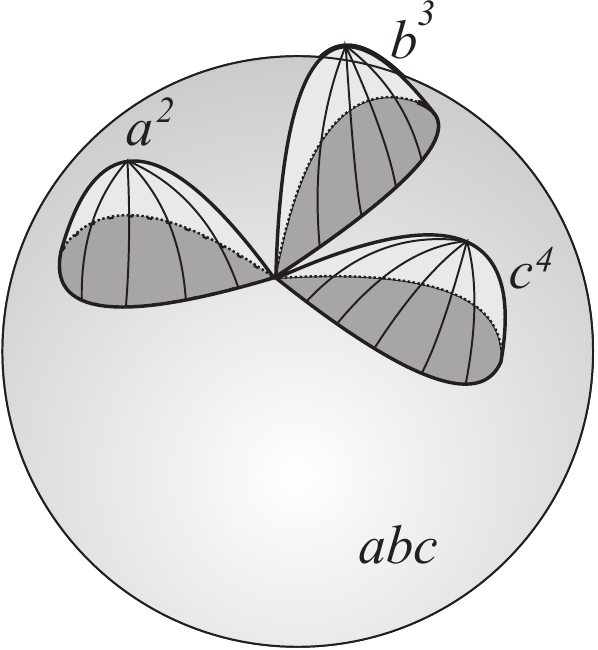}
\end{minipage}
\begin{minipage}{7.2cm}
\includegraphics[scale=0.24]{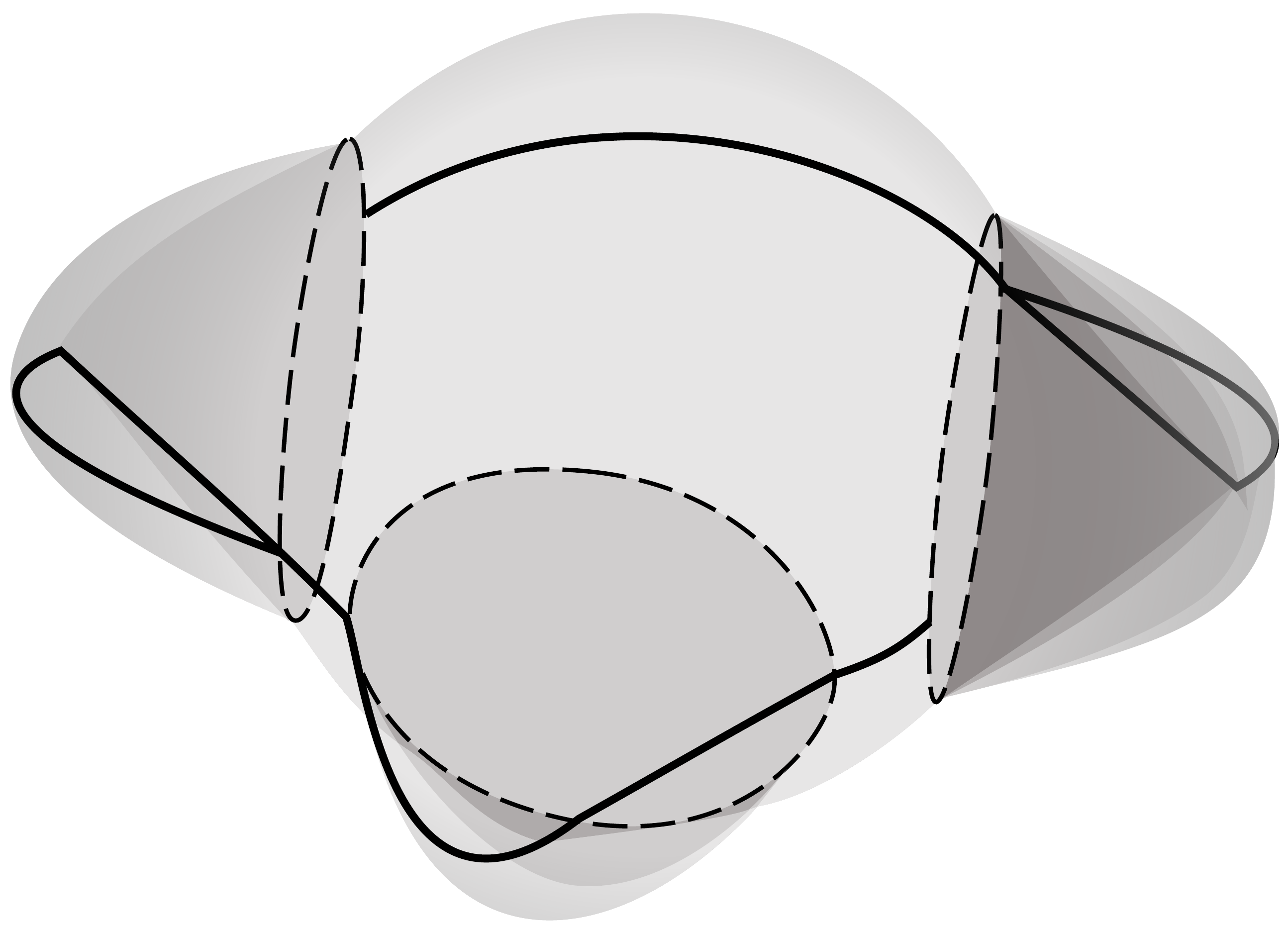}
\end{minipage}
\caption{The space $K_{\PP}$ at the left and the space $X(2,3,4)$ at the right, along with the fixed points of $f$.}\label{figura}
\end{figure}

The space $X=X(l,m,n)$ can be obtained from the surface $S=\{(x,y,z) \,:\, x^2+y^2+z^2$ $=1 \text{ and }x\leq \frac{4}{5}$ and $-\frac{4}{5}\leq y\leq \frac{4}{5} \}$ by attaching three $2$-cells with attaching maps $\phi_a, \phi_b, \phi_c:S^1\to S$ given by $\phi_a(z)=\left(\frac{3}{5}\Re(z^l),-\frac{4}{5},\frac{3}{5}\Im(z^l)\right)$, $\phi_b(z)=\left(\frac{4}{5},\frac{3}{5}\Re(z^m),\frac{3}{5}\Im(z^m)\right)$ and 
$\phi_c(z)=\left(-\frac{3}{5}\Re(z^n),\frac{4}{5},\frac{3}{5}\Im(z^n)\right)$. Let $i_a,i_b, i_c:D^2\to X$ denote the characteristic maps of the cells and let $i_S: S\hookrightarrow X$ be the inclusion.

The maps $f_S:S\to X$, $f_a, f_b, f_c:D^2\to X$ given by $f_S(x,y,z)=i_S(x,y,-z)$, $f_a(w)=i_a(\overline{w})$, $f_b(w)=i_b(\overline{w})$, $f_c(w)=i_c(\overline{w})$, determine a map $f:X\to X$. 

Each connected component of $\Fix(f)$ is homotopy equivalent to $S^1$, as depicted in Figure \ref{figura}. A fixed point class $F$ of $f$ is a union of connected components of $\Fix(f)$. Therefore, $\chi(F)=0$ for every fixed point class $F$ of $f$. An application of Lemma \ref{LemaIndice} yields $i(f,F)=0$. Thus, $N(f)=0$ and we are done.
\end{proof}
\end{teo}

\end{document}